\documentclass[11pt]{amsart}
\makeatletter
\@namedef{subjclassname@2020}{%
  \textup{2020} Mathematics Subject Classification}
\makeatother

\usepackage[a4paper,left=2.4cm,right=2.4cm,top=2.5cm,bottom=2.7cm]{geometry}
\usepackage[T1]{fontenc}
\usepackage[utf8]{inputenc}
\usepackage{lmodern}

\usepackage{amsmath,amssymb,amsthm,mathtools}
\usepackage{microtype}
\usepackage{enumitem}
\usepackage{booktabs}
\usepackage{array}

\usepackage{xcolor}
\usepackage[
    colorlinks=true,
    linkcolor=blue,
    citecolor=blue,
    urlcolor=blue
]{hyperref}

\newtheorem{principle}{Principle}

\theoremstyle{definition}
\newtheorem{definition}{Definition}
\newtheorem{question}{Question}

\theoremstyle{remark}

\title[The Next Fourier Transform]
{The Next Fourier Transform: \\ Rewiring, representation, and Artificial Mathematical Creativity}

\author{Luca Fanelli}
\address{
Ikerbasque, Basque Foundation for Science, Bilbao, Spain;\
EHU Euskal Herriko Unibertsitatea, Departamento de Matemáticas,
Barrio Sarriena s/n, 48940 Leioa, Spain;\
BCAM -- Basque Center for Applied Mathematics,
Alameda de Mazarredo 14--16, 48009 Bilbao, Spain
}
\email{luca.fanelli@ehu.eus}
\email{lfanelli@bcamath.org}
\thanks{ORCID: \href{https://orcid.org/0000-0003-1714-1611}
{0000-0003-1714-1611}.}



\begin{document}

\begin{abstract}
What distinguishes the solution of a difficult mathematical problem from a conceptual revolution? Some mathematical innovations do more than establish new results: they reorganize the relations among problems. We propose to describe such events through changes in the \emph{transfer structure} of mathematics, namely the pattern determining which problems can naturally inform one another, which methods can move between them, and which new questions become accessible.

\noindent The eighteenth-century controversy over the vibrating string provides our central historical case. D'Alembert found a famous model for wave propagation, Bernoulli proposed a modal interpretation, Euler broadened the class of admissible profiles, and Lagrange constructed a finite-dimensional strategy based on modal decomposition. Later, Fourier analysis gave the trigonometric language in which all the previous descriptions can interact.

\noindent  This motivates our \emph{Rewiring Thesis}: a conceptual revolution is a persistent, large-scale reorganization of transfer relations among mathematical problems. The productive representation is one mechanism capable of producing such rewiring, but some stress tests (involving Galois, Riemann, Lebesgue--Schwartz and Wiles) show that it is not the only one. The musical notation provides an independent case in which the representation becomes operationally generative without being informationally complete.

\noindent The resulting framework suggests a stronger criterion for artificial mathematical creativity than  problem solving alone. The relevant question is not merely {\it whether AI can solve a famous conjecture}, but {\it whether it can recognize that an inherited conceptual organization is itself an obstruction and produce a new one that reorganizes several problems at once}. In this sense, the decisive achievement may be the invention of {\bf the next Fourier transform}.
\end{abstract}

\maketitle

\section{Introduction}

Imagine a guitarist, plucking a perfectly flexible string at one point. In the simplest idealization, its initial shape is triangular: if the string has length $\ell$ and is pulled at $x=a$ to height $h$, then

\[
f(x)=
\begin{cases}
\dfrac{h}{a}x, & 0\leq x\leq a,\\[6pt]
\dfrac{h}{\ell-a}(\ell-x), & a\leq x\leq \ell.
\end{cases}
\]

The function is continuous but not differentiable at $x=a$. The normal modes of the string, by contrast, are the smooth functions

\[
\sin\left(\frac{n\pi x}{\ell}\right),
\qquad n=1,2,\ldots,
\]

and today we write, under the appropriate interpretation,

\[
f(x)
=
\sum_{n=1}^{\infty}
b_n
\sin\left(\frac{n\pi x}{\ell}\right).
\]

Every finite partial sum is smooth, while the limiting function has a corner. For us there is no paradox: regularity does not necessarily survive an infinite limiting process, and the behavior of the Fourier coefficients carries information about the regularity of the represented function. In the middle of the eighteenth century, however, the mathematical language in which these facts could comfortably coexist had not yet stabilized.

The well known vibrating-string controversy makes this fragmentation unusually visible. D'Alembert obtained the wave equation and a travelling-wave description; Daniel Bernoulli argued for a superposition of harmonic modes; Euler insisted on a broader notion of admissible function; and Lagrange approached the continuous string through finite systems of interacting masses. The real history is in fact considerably more complicated than a linear succession of discoveries, and we will return to it in Section~\ref{sec:string} (see also \cite{Darrigol2007,Herreman2013}). What matters for the present argument is that several powerful descriptions already existed without yet belonging to a single stable operational language.

The work by Fourier on trigonometric series and heat \cite{Fourier1822} did not suddenly resolve every difficulty left by this controversy; in addition, Fourier did not invent normal modes. Questions of convergence, admissibility and rigor continued to develop afterwards, in particular thanks to Dirichlet's work \cite{Dirichlet1829}. The more interesting historical fact for our purposes is structural: trigonometric representation increasingly became a systematic language in which general data, spectral coefficients, differential operators and different evolution equations could be related to each other through common operations.

This suggests a way of distinguishing two forms of mathematical progress. Solving a problem changes what is known about that problem. A deeper conceptual transformation may instead change how a collection of problems is organized: methods that previously belonged to separate contexts become transferable, distinctions that appeared fundamental become secondary, and new questions become natural.

To formulate this idea, let $\mathcal P_L$ denote a family of mathematical problems available within a conceptual language $L$, where ``language'' includes not only notation but also the objects, representations, transformations, invariants and standard methods through which those problems are approached. We will describe their organization by a transfer structure

\[
\mathfrak T_L=(\mathcal P_L,c_L),
\]

where $c_L(P\to Q)$ represents, schematically, the conceptual cost of transferring an effective way of thinking from $P$ to $Q$. This is not intended as a canonical metric on mathematics. Section~2 develops only the amount of formal structure required to distinguish local progress from changes in the organization of a problem family.

Our central proposal is the following.

\begin{principle}[Rewiring Thesis]\label{pr:rewiring}
A conceptual revolution in mathematics is a persistent, large-scale reorganization of the transfer relations among mathematical problems. It changes which problems can naturally inform one another, which methods can move between them, and which new problems can be naturally formulated.
\end{principle}

The distinction is merely taxonomic and does not pretend to be evaluative. A theorem may be extraordinarily difficult and important while leaving much of the surrounding transfer structure unchanged; conversely, an innovation arising from a relatively specific problem may reorganize an entire domain. In the language developed in later sections,

\[
\text{problem solving changes status},
\qquad
\text{conceptual revolution changes structure}.
\]

This broader distinction between theorem proving and mathematical progress
has an important precedent in Thurston's emphasis on advancing mathematical
understanding rather than identifying progress with the production of formal
results alone \cite{Thurston1994}.
The emphasis on transfer is not itself new. Fisseni, Sarikaya and
Schröder have analysed mathematical innovation through the transfer of
tools, structural concepts and representations between apparently remote
fields \cite{FisseniSarikayaSchroeder2023}. Likewise, the distinction
between exploring a given conceptual space and transforming the space
itself has a substantial precedent in Boden's account of exploratory and
transformational creativity \cite{Boden1998}. More recently, Jutla and
Sharma have connected the history of the vibrating string and Fourier
analysis with the problem of artificial mathematical creativity
\cite{JutlaSharma2026}. The present proposal differs from these accounts
in the object it attempts to make explicit: a directed 
transfer structure between problems, together with separate notions of the amplitude,
breadth and persistence of its deformation, and of expansion of the
accessible problem space. The aim is not merely to observe that concepts,
representations or methods can move between domains, but to ask when such
changes reorganize the relations among several mathematical problems
strongly and persistently enough to count as a change in mathematical
structure rather than only as progress within it.

Fourier analysis also suggests one especially important mechanism for producing such a transformation. A new representation may matter not because it preserves more information, but because it changes which operations are accessible.

\begin{principle}[Productive Representation Principle]\label{pr}
A representation is conceptually productive when it does more than encode an existing object. It changes the operations that can naturally be performed on that object and, through those operations, may change the transfer structure of the surrounding domain.
\end{principle}

The logical relation between the two principles is important. The Rewiring Thesis is the central claim of the paper. Productive representation is one of the mechanisms capable of producing rewiring; it is not our definition of conceptual revolution. In the sequel, some {\it stress tests} will show that new structural mediators, new mathematical objects and enlargements of the admissible universe can produce comparable reorganizations without being naturally described as changes of representation.

This is the point in which Music enters the argument. \textbf{We do not claim that mathematics and music undergo the same kinds of conceptual revolution, nor that a musical score is analogous to a Fourier transform}. Musical notation provides a more specific comparison. Within the Western tradition, notation gradually transformed aspects of an essentially temporal practice into symbolic objects that could be inspected, compared and manipulated on a surface. Its development interacted with compositional practice rather than merely recording it \cite{BusseBerger2002,Grier2021}. At the same time, a score is clearly incomplete as a representation of the performance.

Music therefore provides an independent test of the Productive Representation Principle. It shows especially clearly that operational power and informational completeness are different properties: a representation may become generative precisely because it stabilizes some distinctions while leaving others aside. Recent philosophical work has similarly emphasized the productive role of notation rather than treating it as a neutral record \cite{ColyvanHall2026}.

This selectivity leads to a complementary mathematical question. If

\[
R:X\longrightarrow Z
\]

is a representation, which operations and tasks survive passage through $R$, and which do not? Section~\ref{sec:blindspots} will show that these two possibilities can be formulated through factorization: productive operations factor effectively through a representation, whereas a task that distinguishes elements collapsed by the same representation cannot be recovered exactly from it. Affordances and blind spots are therefore two sides of the same representational choice.

These issues become particularly relevant for artificial intelligence. Modern systems already learn internal representations, so the interesting contrast cannot be between humans who invent representations and machines that merely manipulate fixed ones. We will instead distinguish

\[
\text{representation learning}
\longrightarrow
\text{representation invention}
\longrightarrow
\text{representational revolution}.
\]

Learning an effective encoding may improve performance without changing mathematics. Representation invention requires a new externalizable system of objects or relations. Representational revolution requires still more: the new structure must reorganize relations among mathematical problems.

This leads to a concrete question about artificial conceptual creativity.

\begin{question}[Representational Escape Problem]\label{q:escape}
Can an artificial system recognize that the representation in which a family of problems is currently posed is itself an obstruction, replace it autonomously by another representation, and consequently reorganize several previously distinct problems at once?
\end{question}

Representational escape is deliberately stronger than improved search inside a fixed framework, but it is also narrower than conceptual revolution in general. A successful escape becomes revolutionary only if the resulting structure produces the broad and persistent reorganization described by the Rewiring Thesis. Conversely, the stress tests developed later in the paper show that conceptual rewiring could occur through mechanisms other than representational escape.

This distinction sharpens a question already present in recent discussions of AI and mathematical practice. Mathematical activity includes much more than proof production: theory building, communication, verification, digestion and canonicalization also matter \cite{CommelinEtAl2026, Tao2026}. Our concern is one specific part of that broader landscape. We ask what observable mathematical consequences would justify saying that an artificial system had changed the conceptual organization of a domain rather than merely become more effective within it.

The paper is organized as follows. Section~2 introduces transfer structures and a parameterized notion of conceptual revolution based on amplitude, breadth and persistence, with expansion as an additional feature. Section~\ref{sec:string} reconstructs the fragmented pre-Fourier situation of the vibrating string, and Section~\ref{sec:fourier-rewiring} examines Fourier analysis as the central positive case. Section~\ref{sec:notation} uses musical notation as an independent test of productive representation, while Section~\ref{sec:blindspots} develops the complementary theory of representational blind spots. Section~\ref{sec:stress-tests} then tests the broader Rewiring Thesis against Galois theory, Riemannian geometry, the enlargements associated with Lebesgue and Schwartz, and the modular architecture surrounding Fermat's Last Theorem. Finally, Section~\ref{sec:ai} returns to artificial intelligence and develops representational escape as one operational test of conceptual creativity.

The hierarchy of the argument is therefore simple: the Rewiring Thesis is the central proposal; productive representation is one mechanism capable of producing rewiring; and representational escape asks whether an artificial system could autonomously produce such a transformation from within an inherited conceptual framework.

%%%%%%%%%%%%%%%%%%%%%%%%%%%%%%%%%%%%%%%%%%%%%%%%%%%%
%%%%%%%%%%%%%%%%%%%%%%%%%%%%%%%%%%%%%%%%%%%%%%%%%%%%
%%%%%%%%%%%%%%%%%%%%%%%%%%%%%%%%%%%%%%%%%%%%%%%%%%%%
%%%%%%%%%%%%%%%%%%%%%%%%%%%%%%%%%%%%%%%%%%%%%%%%%%%%

\section{Transfer structure and conceptual revolution}
\label{sec:framework}

The expression ``conceptual revolution'' is easy to use and difficult to make precise. A theorem, object or technique may be extremely important without reorganizing a field, while a development that initially appears local may eventually alter how many different problems interact. Our proposal is therefore to focus not only on mathematical results, but on the relations among problems.

\subsection{Problems, languages and transfer}

For a conceptual language \(L\), let

\[
\mathcal P_L=\{P_1,P_2,\ldots\}
\]

denote a family of mathematical problems available to a mathematical community. The word ``problem'' is used in a broad sense: it may refer to a specific conjecture, the study of a class of equations, a classification problem, the construction of examples, or a recurring mathematical task. By a language \(L\) we mean more than notation. It includes the objects regarded as natural, the representations used for them, the standard transformations and invariants, the available methods of proof, and the analogies mathematicians know how to exploit.

For \(P,Q\in\mathcal P_L\), let us introduce a non-negative function

\[
c_L:
\mathcal P_L\times\mathcal P_L
\longrightarrow
[0,+\infty)
\]

and write

\[
c_L(P\to Q):=c_L(P,Q).
\]

We call \(c_L(P\to Q)\) the \emph{transfer cost from \(P\) to \(Q\) in the language \(L\)}. It represents, schematically, the difficulty of transferring an effective way of thinking from \(P\) to \(Q\). If a method developed for \(P\) adapts naturally to \(Q\), the cost should be small; if progress on \(P\) provides little usable structure for \(Q\), it should be large.

We assume

\[
c_L(P,P)=0
\qquad
\text{for every }P\in\mathcal P_L.
\]

Since transfer need not be symmetric, we do not assume here that

\[
c_L(P\to Q)=c_L(Q\to P),
\]

nor do we impose a triangle inequality. Thus \(c_L\) does not need to be understood as a metric, and later references to problems becoming ``closer'' or changing their ``distance'' are just metaphors for changes in directed transfer cost.

We also do not assume that transfer costs possess a canonical and absolute numerical scale. Whenever two languages are compared, their cost functions are expressed on a common comparison scale, appropriate to that comparison. A common positive rescaling changes the numerical thresholds below but not the structural classification, whenever those thresholds are rescaled accordingly.

\begin{definition}[Transfer structure]
Given a conceptual language \(L\), the corresponding \emph{transfer structure} is

\[
\mathfrak T_L=(\mathcal P_L,c_L).
\]

\end{definition}

The dependence on \(L\) is essential. A change of language may alter both the costs of transferring methods among existing problems and the collection of problems that can naturally be formulated.

\subsection{Innovation and transfer deformation}

\begin{definition}[Innovation]
An \emph{innovation} \(\mathcal I\) is a transformation

\[
\mathcal I:
\mathfrak T_L
\longmapsto
\mathfrak T_{L_{\mathcal I}},
\]

that is,

\[
\mathcal I:
(\mathcal P_L,c_L)
\longmapsto
(\mathcal P_{L_{\mathcal I}},c_{L_{\mathcal I}}),
\]

together with an injective tracking map

\[
\iota_{\mathcal I}:
\mathcal P_L
\longrightarrow
\mathcal P_{L_{\mathcal I}}
\]

that identifies the continuation of each old problem inside the new transfer structure.
\end{definition}

Injectivity is important. A new language may reveal that two previously distant problems are instances of a common mathematical structure, but this cannot make the original tasks literally identical. Such unification should be reflected by a reduction in their transfer costs, not by collapsing them into a single vertex.

Thus, if \(P\neq Q\), then

\[
\iota_{\mathcal I}(P)
\neq
\iota_{\mathcal I}(Q),
\]

even when the innovation makes the relation between them conceptually immediate.

For a fixed comparison between \(L\) and \(L_{\mathcal I}\), the two cost functions are  understood to be on the common comparison scale mentioned above.

\begin{definition}[Transfer deformation function]\label{def:transfer}
Given an innovation \(\mathcal I\), define

\[
\Delta_{\mathcal I}(P,Q)
:=
c_{L_{\mathcal I}}
\bigl(
\iota_{\mathcal I}(P),
\iota_{\mathcal I}(Q)
\bigr)
-
c_L(P,Q),
\qquad
P,Q\in\mathcal P_L.
\]

\end{definition}

Thus

\[
\Delta_{\mathcal I}(P,Q)<0
\]

means that the transfer from \(P\) to \(Q\) has become easier, while

\[
\Delta_{\mathcal I}(P,Q)>0
\]

means that it has become less natural in the new language. We do not expect \(\Delta_{\mathcal I}\leq0\) for every pair. An innovation may reveal some relations while obscuring others. The appropriate notion is therefore a \emph{deformation}, not necessarily a contraction, of transfer structure.

An innovation may also enlarge the problem space. Define

\[
\mathcal N_{\mathcal I}
=
\mathcal P_{L_{\mathcal I}}
\setminus
\iota_{\mathcal I}(\mathcal P_L).
\]

Because every old problem has a distinct tracked continuation, the image
\(\iota_{\mathcal I}(\mathcal P_L)\) records the old problem family as represented
in the new language. We adopt the convention that alternative reformulations
of an old problem are not counted as distinct new problems unless they define
a genuinely new mathematical task. Under this convention, the elements of
\(\mathcal N_{\mathcal I}\) are precisely those problems in the new language
that have no antecedent among the problems available in \(L\). Accordingly,

\[
\mathcal N_{\mathcal I}\neq\varnothing
\]

expresses genuine expansion of the naturally available problem space rather
than an artifact of identifying or merely reformulating previously existing
problems.

\subsection{Substantiality, non-locality and persistence}

Let

\[
F=\{P_1,\ldots,P_N\}\subset\mathcal P_L,
\qquad
N\geq2,
\]

be a finite reference family. The parameters introduced below determine the resolution at which an innovation is examined and are fixed for the comparison rather than chosen a posteriori to make a candidate innovation satisfy the criteria.

The reference family \(F\) is itself part of that scale. It should be specified independently of the particular deformation one intends to exhibit; otherwise breadth could be inflated by selecting only problems already known to be strongly affected. In historical applications, we therefore prefer families defined by the mathematical situation before the candidate innovation is analysed, or by an independently stated criterion of scope.

For \(\varepsilon>0\), define

\[
A_{\mathcal I}(F,\varepsilon)
=
\left\{
(P,Q)\in F\times F:
P\neq Q,\ 
|\Delta_{\mathcal I}(P,Q)|\geq\varepsilon
\right\}.
\]

\begin{definition}[Substantial deformation]
The innovation \(\mathcal I\) is \emph{\(\varepsilon\)-substantial on \(F\)} if

\[
A_{\mathcal I}(F,\varepsilon)\neq\varnothing.
\]

\end{definition}

Thus \(\varepsilon\) fixes the amplitude at which a transfer deformation counts as significant.

Define the corresponding set of affected problems by

\[
V_{\mathcal I}(F,\varepsilon)
=
\left\{
P\in F:
\exists Q\in F
\text{ such that }
(P,Q)\in A_{\mathcal I}(F,\varepsilon)
\text{ or }
(Q,P)\in A_{\mathcal I}(F,\varepsilon)
\right\}.
\]

\begin{definition}[Non-local deformation]
For \(0<\delta\leq1\), the innovation \(\mathcal I\) is
\emph{\((\varepsilon,\delta)\)-non-local on \(F\)} if

\[
\frac{|V_{\mathcal I}(F,\varepsilon)|}{|F|}
\geq\delta.
\]

\end{definition}

The parameter \(\delta\) measures breadth. Since \(\delta>0\),

\[
(\varepsilon,\delta)\text{-non-locality}
\quad\Longrightarrow\quad
\varepsilon\text{-substantiality}.
\]

Amplitude and breadth are therefore conceptually distinct but not logically independent conditions.

To formulate persistence, suppose that the innovation is subsequently incorporated into languages

\[
L_{\mathcal I}=L_0,L_1,L_2,\ldots
\]

with injective tracking maps

\[
\iota_n:
\mathcal P_L
\longrightarrow
\mathcal P_{L_n},
\qquad
\iota_0=\iota_{\mathcal I}.
\]

Assume that the corresponding transfer costs are expressed on a common comparison scale and define

\[
\Delta_{\mathcal I}^{(n)}(P,Q)
=
c_{L_n}
\bigl(
\iota_n(P),\iota_n(Q)
\bigr)
-
c_L(P,Q).
\]

Then

\[
\Delta_{\mathcal I}^{(0)}
=
\Delta_{\mathcal I}.
\]

\begin{definition}[Persistent deformation at scale \(\varepsilon\)]
Let

\[
\varnothing\neq B\subset F\times F.
\]

The deformation on \(B\) is \emph{\(\varepsilon\)-persistent} if there exists \(n_0\) such that

\[
|\Delta_{\mathcal I}^{(n)}(P,Q)|
\geq
\varepsilon
\]

for every \((P,Q)\in B\) and every \(n\geq n_0\).
\end{definition}

Persistence therefore means more than remaining detectably different from the original transfer structure. The deformation must remain significant at the same comparison scale \(\varepsilon\) at which it initially counted as substantial. A large initial reorganization that eventually decays below that scale does not qualify as persistent at that scale.

For \(B\subset F\times F\), define its incident set

\[
V(B)
=
\left\{
P\in F:
\exists Q\in F
\text{ such that }
(P,Q)\in B
\text{ or }
(Q,P)\in B
\right\}.
\]

We can now state the central definition.

\begin{definition}[Conceptual revolution]
Fix a finite reference family

\[
F\subset\mathcal P_L,
\qquad
|F|\geq2,
\]

a significance threshold \(\varepsilon>0\), and breadth parameters

\[
0<\delta\leq1,
\qquad
0<\delta_*\leq1.
\]

An innovation

\[
\mathcal I:
\mathfrak T_L
\longmapsto
\mathfrak T_{L_{\mathcal I}}
\]

is an
\emph{\((\varepsilon,\delta,\delta_*)\)-conceptual revolution on \(F\)}
if:

\begin{enumerate}

\item \(\mathcal I\) is \((\varepsilon,\delta)\)-non-local on \(F\);

\item there exists a nonempty set

\[
B\subset
A_{\mathcal I}(F,\varepsilon)
\]

such that the deformation on \(B\) is \(\varepsilon\)-persistent and

\[
\frac{|V(B)|}{|F|}
\geq
\delta_*.
\]

\end{enumerate}

In the expansive case,

\[
\mathcal N_{\mathcal I}\neq\varnothing.
\]

\end{definition}

The parameters separate the dimensions relevant to the framework. The threshold \(\varepsilon\) fixes the amplitude at which a deformation counts as significant; \(\delta\) measures the breadth of the initial deformation; and \(\delta_*\) measures the breadth of a core whose deformation remains above the same significance threshold after the assimilation. Expansion records the appearance of naturally accessible problems outside the tracked image of the original problem family.

The definition is deliberately parameterized. The thresholds are not universal constants and may depend on the domain and historical scale under consideration. What matters formally is that the reference family, comparison scale and parameters are fixed independently of the deformation one wishes to establish. Otherwise arbitrarily small thresholds or a selectively chosen family could make the criteria vacuous.

Accordingly, below we use the unqualified expression \emph{conceptual revolution} as qualitative shorthand for this parameterized notion at a contextually meaningful, nontrivial scale. The formalism is intended to distinguish amplitude, breadth, persistence and expansion, not to assign numerical rankings to historical innovations.

\subsection{Problem solving and structural change}

For the present comparison, let

\[
\sigma_L:
\mathcal P_L
\longrightarrow
\{\text{open},\text{solved}\}
\]

describe the status of a problem. Solving \(P\) changes

\[
\sigma_L(P):
\text{open}
\longrightarrow
\text{solved},
\]

but need not substantially alter the surrounding transfer structure. A conceptual revolution acts instead on

\[
\mathfrak T_L=(\mathcal P_L,c_L).
\]

In short,

\[
\text{problem solving changes status},
\qquad
\text{conceptual revolution changes structure}.
\]

This is just a distinction, not a ranking of the two concepts. A theorem may be extraordinarily difficult and important while producing little structural reorganization; an innovation arising from a more local question may eventually alter relations across a wide domain.

\subsection{Productive representations}

Representations provide one important mechanism for such reorganization. Let

\[
R:X\longrightarrow Z
\]

be a representation and \(A:X\to X\) an operation. Suppose there exists an operation \(B:Z\to Z\) such that

\[
R\circ A
=
B\circ R.
\]

This factorization means that the \emph{represented effect} of \(A\) can be obtained from the represented object \(R(x)\): one may pass from \(R(x)\) to \(R(Ax)\) by applying \(B\).

The factorization alone, however, is not sufficient for productivity. A constant representation, for example, can satisfy such an identity trivially while preserving no structure relevant to the mathematical tasks under consideration. For a factorization to constitute a productive affordance, the represented operation must preserve distinctions which are relevant to those tasks and yield a genuine operational advantage (for example, by making a family of calculations, comparisons or transfers substantially easier).

Thus the diagram

\[
R\circ A=B\circ R
\]

provides a structural mechanism through which productivity may occur; productivity itself is relative to the mathematical operations and tasks that the representation makes accessible.

If the induced operation \(B\) supplies such an advantage across a sufficiently broad family of problems, then \(R\) can deform the surrounding transfer structure. This gives the formal connection between operational representation and the Productive Representation Principle, Principle~\ref{pr}.

Expansion through representation is expressed by

\[
\iota_{\mathcal I}(\mathcal P_L)
\subsetneq
\mathcal P_{L_{\mathcal I}},
\]

equivalently,

\[
\mathcal N_{\mathcal I}\neq\varnothing.
\]

A productive representation does not necessarily improve every transfer. Some operations may become easier while others become harder or less visible. This possibility is already encoded by the fact that

\[
\Delta_{\mathcal I}(P,Q)
\]

can have either sign.

The same factorization viewpoint will clarify one class of representational limitations in Section~\ref{sec:blindspots}: a task that does not factor through \(R\) cannot be recovered exactly from \(R(x)\) alone. Fiber structure therefore characterizes one form of information loss. Operational productivity, by contrast, depends not only on which distinctions a representation preserves, but also on how the preserved structure organizes the operations relevant to the task.

The roles of representation and reformulation have been studied from several perspectives (see e.g. \cite{Hunt2025,Ippoliti2022,Schlimm2025}). Our specific proposal is to connect representational productivity with changes in the transfer structure.

\subsection{Scope of the framework}

Several limitations should remain explicit. Transfer costs are not assumed to possess a canonical absolute numerical scale; the formalism requires only a comparison scale for the languages under consideration. Moreover, not every form of mathematical progress is a rewiring: solving problems, proving estimates, constructing examples, developing techniques, classifying objects and making theories rigorous remain genuine forms of progress without necessarily producing a conceptual revolution.

Conceptual revolutions may also be historically distributed across many people and decades, and one representation not necessarily permanently replaces another. Most importantly, not every revolution is representational. New invariants, structural objects, symmetry principles, abstractions or enlarged classes of admissible objects may also reorganize transfer relations. These mechanisms need not be mutually exclusive: the same historical development may admit more than one useful description. Productive representation is one mechanism for rewiring, not its definition.

The framework is related to broader accounts of mathematical and conceptual change. Kuhn provides an important background notion of scientific revolution in \cite{Kuhn1970}, while Lakatos analysed the joint evolution of definitions, counterexamples and proof strategies within mathematical practice in \cite{Lakatos1976}. Geometrical accounts of conceptual change have been developed through conceptual spaces (see \cite{Gardenfors2000,GardenforsZenker2013}) and applied to mathematical examples in \cite{PetersenZenker2014}. Our object is narrower: the vertices are mathematical problems, and the structure of interest is the directed transfer of methods, representations and forms of reasoning among them.

For the same reason, \(\mathfrak T_L\) should not be confused with the classical notion of a \emph{problem space} in cognitive science \cite{NewellSimon1972}, which describes states and possible moves within the solution of a problem. A transfer structure concerns relations between different mathematical problems.

The vibrating string provides an unusually clean first test. By the middle of the eighteenth century, the equation, explicit descriptions of the motion and the idea of normal-mode decomposition were already available, yet the passages among these descriptions had not stabilized into a common operational language. Let us now turn to that fascinating history.

%%%%%%%%%%%%%%%%%%%%%%%%%%%%%%%%%%%%%%%%%%%%%%%%%
%%%%%%%%%%%%%%%%%%%%%%%%%%%%%%%%%%%%%%%%%%%%%%%%%
%%%%%%%%%%%%%%%%%%%%%%%%%%%%%%%%%%%%%%%%%%%%%%%%%

\section{The vibrating string before Fourier: one problem, several languages}
\label{sec:string}

The vibrating-string controversy is one of the most extensively studied episodes in the history of analysis. Its development has been reconstructed in detail by Truesdell, Darrigol, Guilbaud and Jouve, Herreman, and Jouve \cite{Truesdell1960,Darrigol2007,GuilbaudJouve2009,Herreman2013,Jouve2017}. Herreman has examined the conditions under which trigonometric series acquired a new representational status, while Darrigol has emphasized the acoustic origins of harmonic analysis and the complex continuity between the work of Bernoulli, Euler, Lagrange and Fourier.

Our purpose is not to provide another history of the controversy, nor to suggest that eighteenth-century mathematicians were thinking in the terminology introduced in Section~2. We use the episode instead as a historical test of that terminology. The central feature of the pre-Fourier situation, from this viewpoint, is not a lack of powerful mathematics. Important solutions, representations and limiting procedures were already available. What remained unstable were the passages between them. The mathematics was not empty; it was fragmented.

\subsection{A small problem family}

Consider the finite family

\[
F_{\mathrm{string}}
=
\left\{
P_{\mathrm{prop}},
P_{\mathrm{mode}},
P_{\mathrm{repr}},
P_{\mathrm{coeff}},
P_{\mathrm{lim}}
\right\},
\]

where

\[
\begin{aligned}
P_{\mathrm{prop}}
&:
\text{determine the motion of the string from its initial data},\\[3pt]
P_{\mathrm{mode}}
&:
\text{describe the motion in terms of normal modes},\\[3pt]
P_{\mathrm{repr}}
&:
\text{represent a sufficiently general initial profile},\\[3pt]
P_{\mathrm{coeff}}
&:
\text{recover the modal amplitudes from that profile},\\[3pt]
P_{\mathrm{lim}}
&:
\text{pass from finite mechanical systems to the continuous string}.
\end{aligned}
\]

These are modern labels for historical mathematical difficulties. We do not claim that d'Alembert, Euler, Bernoulli or Lagrange formulated them in precisely this way. The point is to keep the same reference family fixed while asking which transfers different mathematical languages made available.

As in Section~2, we assign no numerical values to transfer costs. Saying that $c_L(P\to Q)$ is large means only that the available language $L$ provides no stable general procedure for transferring a solution, representation or method from $P$ to $Q$ without substantial additional conceptual work.

\subsection{D'Alembert: propagation becomes accessible}

D'Alembert's work led to the equation now written as

\[
u_{tt}=c^2u_{xx},
\]

together with the travelling-wave representation \cite{dAlembert1749a,dAlembert1749b}

\[
u(x,t)
=
\Phi(x-ct)+\Psi(x+ct).
\]

For zero initial velocity, in modern notation,

\[
u(x,t)
=
\frac12
\bigl(
f(x-ct)+f(x+ct)
\bigr),
\]

with the appropriate reflected extension for fixed endpoints. Once an admissible initial profile $f$ is given, propagation can therefore be described directly through translations of that profile.

This makes $P_{\mathrm{prop}}$ comparatively accessible, but it does not automatically stabilize the other transfers in $F_{\mathrm{string}}$. The travelling-wave description does not determine modal coefficients or explain why the same motion should admit a decomposition into harmonic modes. Moreover, the admissible class of profiles $f$ was itself part of the controversy \cite{Jouve2017}.

Thus d'Alembert provides an early example of the distinction central to this paper: \textbf{a problem may acquire a powerful solution before its different descriptions become mutually transferable}.

\subsection{Bernoulli: a visible relation before a stable transfer}

Daniel Bernoulli approached the same physical system through its elementary modes

\[
\sin\left(\frac{n\pi x}{\ell}\right),
\qquad
n=1,2,\ldots,
\]

and argued that a general vibration should be understood as a superposition of such harmonic motions \cite{Bernoulli1755a,Bernoulli1755b}. In modern notation, for zero initial velocity,

\[
u(x,t)
=
\sum_{n=1}^{\infty}
b_n
\sin\left(\frac{n\pi x}{\ell}\right)
\cos\left(\frac{n\pi ct}{\ell}\right),
\]

so that formally at $t=0$,

\[
f(x)
=
\sum_{n=1}^{\infty}
b_n
\sin\left(\frac{n\pi x}{\ell}\right).
\]

In terms of $F_{\mathrm{string}}$, Bernoulli was proposing a general connection between $P_{\mathrm{repr}}$ and $P_{\mathrm{mode}}$. But recognizing that such a connection should exist is not the same as possessing a stable mathematical transfer. For a general profile $f$, one still had to determine whether the representation existed, how the coefficients $b_n$ should be obtained, in what sense the equality held, and under what conditions the series converged to the prescribed profile.

Consequently, the transfer

\[
P_{\mathrm{repr}}
\longrightarrow
P_{\mathrm{coeff}}
\]

was not yet available as a general systematic procedure, while the passage

\[
P_{\mathrm{coeff}}
\longrightarrow
P_{\mathrm{mode}}
\]

was available only at a formal level: what remained unsettled was the general
validity of the modal reconstruction, including the existence, convergence and
interpretation of the resulting infinite series.

The distinction is conceptually important: a mathematical relation may be \emph{visible} before it becomes \emph{transferable}. Bernoulli's physical interpretation identified a connection that the available analytical language still had not stabilized.

\subsection{Euler: changing the admissible problem space}

Euler's intervention reveals a different mechanism. A physical string may be assigned initial shapes that do not belong naturally to a narrow class of curves represented by a single analytic expression. Euler insisted on admitting substantially more general curves into the discussion \cite{Darrigol2007,Euler1750, Jouve2017}. The controversy therefore concerned not only how to solve the wave equation, but also what kind of object an admissible function could be.

In the language of Section~2, this affects not only transfer relations but the collection of naturally legitimate problems,

\[
\mathcal P_L
\longrightarrow
\mathcal P_{L_{\mathcal I}}.
\]

This schematic notation should not be read as a literal set-theoretic reconstruction of eighteenth-century mathematics. It expresses a structural point: a language allowing only a restricted class of curves supports fewer initial-value problems than one admitting piecewise-defined or otherwise irregular profiles.

Euler's contribution therefore illustrates an expansive mechanism before it provides a universal method of solution: changing the admissible objects changes which problems can naturally be posed.

\subsection{Lagrange: a finite route toward the continuum}

Lagrange approached the continuous string through finite systems of interacting masses and then considered the passage to the continuum \cite{Darrigol2007, Lagrange1759}. In simplified modern notation, let

\[
d_N=\frac{\ell}{N+1}
\]

denote the distance between consecutive masses. Such a system can be written as

\[
\ddot y_k
=
\alpha_N
\bigl(
y_{k-1}-2y_k+y_{k+1}
\bigr),
\qquad
k=1,\ldots,N,
\]

with suitable boundary conditions. For the discrete system to approach the wave equation in the continuum limit, the coupling must scale with the lattice spacing so that

\[
\alpha_N d_N^2 \longrightarrow c^2
\qquad
\text{as }N\longrightarrow\infty.
\]

Indeed,

\[
\frac{y_{k-1}-2y_k+y_{k+1}}{d_N^2}
\]

is the standard finite-difference approximation of the second spatial derivative. For each fixed \(N\), the system is finite-dimensional and can be diagonalized and decomposed into finitely many normal modes.

The difficult passage is

\[
N\longrightarrow\infty.
\]

For finite $N$, only finitely many modes occur; in the continuous limit, one seeks to represent a much larger class of initial profiles by infinitely many modes. Lagrange's construction therefore creates an important route

\[
P_{\mathrm{lim}}
\longrightarrow
P_{\mathrm{mode}}
\longrightarrow
P_{\mathrm{repr}},
\]

but the last passage remains unstable. The finite-dimensional diagonalization strongly suggests the infinite modal description without by itself supplying a general theory explaining when the limiting series represents the prescribed profile.

The transition from finite modal decomposition to infinite representation is therefore one of the clearest places where the available transfers had not yet stabilized.

\subsection{The triangular string as a transfer test}

The remaining fragmentation becomes especially transparent for the triangular profile introduced in the Introduction:

\[
f(x)=
\begin{cases}
\dfrac{h}{a}x,
&0\leq x\leq a,\\[6pt]
\dfrac{h}{\ell-a}(\ell-x),
&a\leq x\leq \ell.
\end{cases}
\]

It is continuous but has a corner at $x=a$,

\[
f'(a^-)=\frac{h}{a},
\qquad
f'(a^+)=-\frac{h}{\ell-a}.
\]

From the modern point of view, the same profile can both be propagated by the wave equation and represented by normal modes,

\[
f(x)
=
\sum_{n=1}^{\infty}
b_n
\sin\left(\frac{n\pi x}{\ell}\right),
\]

with

\[
b_n
=
\frac{2h\ell^2}
{a(\ell-a)n^2\pi^2}
\sin\left(\frac{n\pi a}{\ell}\right).
\]

In the symmetric case $a=\ell/2$, we have

\[
b_n
=
\frac{8h}{\pi^2n^2}
\sin\left(\frac{n\pi}{2}\right).
\]

Every individual mode, and every finite partial sum

\[
S_N(x)
=
\sum_{n=1}^{N}
b_n
\sin\left(\frac{n\pi x}{\ell}\right),
\]

is smooth, whereas the limiting profile has a corner. In the symmetric case the nonzero coefficients satisfy

\[
|b_n|
=
\frac{8h}{\pi^2n^2},
\qquad
n\ \text{odd},
\]

and hence decay like $n^{-2}$.

Today the relation between spectral decay and regularity belongs to the ordinary language of analysis. It should not be projected backwards onto the eighteenth-century controversy. For our purposes, however, it illustrates something that becomes possible once Fourier coefficients are stable mathematical objects: a new question emerges naturally,

\[
P_{\mathrm{reg}}
=
\text{determine regularity from spectral decay}.
\]

This question is not naturally present in the original physical formulation of the string. A geometric irregularity in physical space has become information distributed through the asymptotic behavior of spectral coordinates. In other words, this example anticipates both features central to the next section: existing problems become linked through new transfers, and the new language makes additional questions natural.

\subsection{A fragmented transfer structure}

The above contributions should not be forced into a linear sequence or treated as competing solutions of a single problem. They act on different parts of the same reference family. Schematically,

\[
\begin{aligned}
\mathcal I_{\mathrm{dA}}
&:
P_{\mathrm{prop}}
\quad\text{becomes accessible},\\[3pt]
\mathcal I_{\mathrm{B}}
&:
P_{\mathrm{repr}}
\longleftrightarrow
P_{\mathrm{mode}}
\quad\text{is proposed as a general connection},\\[3pt]
\mathcal I_{\mathrm{E}}
&:
\mathcal P_L
\quad\text{is enlarged through a broader notion of admissible profile},\\[3pt]
\mathcal I_{\mathrm{La}}
&:
P_{\mathrm{lim}}
\longrightarrow
P_{\mathrm{mode}}
\quad\text{is made substantially easier}.
\end{aligned}
\]

D'Alembert made propagation accessible; Bernoulli identified a general modal connection; Euler enlarged the admissible class of objects; and Lagrange supplied a controlled finite-dimensional route toward modal decomposition and the continuum. What was still missing was neither another isolated solution, nor  the discovery of normal modes or the mere suggestion of an infinite trigonometric expansion: the missing ingredient was a sufficiently general and operational language in which representation, coefficient extraction, modal decomposition and evolution could become systematically transferable. The pre-Fourier situation was therefore mathematically rich but structurally fragmented.

The next section asks whether Fourier analysis supplied the broad and persistent reorganization needed to turn that fragmentation into a common mathematical architecture.

%%%%%%%%%%%%%%%%%%%%%%%%%%%%%%%%%%%%%%%%%%%%%%%
%%%%%%%%%%%%%%%%%%%%%%%%%%%%%%%%%%%%%%%%%%%%%%%
%%%%%%%%%%%%%%%%%%%%%%%%%%%%%%%%%%%%%%%%%%%%%%%

\section{Fourier and the rewiring of the problem space}
\label{sec:fourier-rewiring}

We now arrive at the first substantial historical test of the framework introduced in Section~\ref{sec:framework}. The previous section showed that the mathematics preceding Fourier was already rich: d'Alembert had made propagation accessible, Bernoulli had proposed a modal description, Euler had enlarged the class of admissible profiles, and Lagrange had constructed a finite-dimensional route toward modal decomposition and the continuum. What remained comparatively unstable was the operational relation among these descriptions.

Here it is essential to stress the historical caution. Fourier did not invent trigonometric modes, nor the coefficient extraction by trigonometric integration was new with him. Important antecedents occur in the eighteenth-century development of the vibrating-string problem and harmonic analysis \cite{Darrigol2007,Herreman2013}. Moreover, the work by Fourier did not immediately settle questions of convergence and admissibility; these continued to develop afterwards, in particular thanks to Dirichlet's work \cite{Fourier1822,Dirichlet1829}. Our claim is therefore not one of simple priority. The relevant change is the increasing generality, systematic deployment and representational status of trigonometric expansion as an operational language connecting data, coefficients, operators and evolution equations.

\subsection{From modes to coordinates}

Let

\[
\mathcal I_{\mathcal F}:
\mathfrak T_L
\longmapsto
\mathfrak T_{L_{\mathcal F}}
\]

denote the Fourier innovation, where \(L_{\mathcal F}\) is a language in which trigonometric coefficients are treated systematically as coordinates of a function.

This is the conceptual change relevant to our argument. Trigonometric modes and coefficient formulas had important pre-Fourier antecedents. What becomes increasingly stable in the Fourier architecture is their use as a general coordinate system for classes of data.

In modern terms, for \(f\in L^2(0,\ell)\) one has sine coefficients

\[
b_n
=
\frac{2}{\ell}
\int_0^\ell
f(x)
\sin\left(\frac{n\pi x}{\ell}\right)\,dx,
\]

and the sine functions form a complete orthogonal system in \(L^2(0,\ell)\). Thus the coefficients determine the sine-series representation of \(f\) in \(L^2\), with stronger convergence statements available under additional regularity assumptions.

In the notation of Section~\ref{sec:string}, this stabilizes the relation among

\[
P_{\mathrm{repr}},
\qquad
P_{\mathrm{coeff}},
\qquad
P_{\mathrm{mode}}.
\]

The historical claim should not be read as saying that these transfers were absent before Fourier. Rather, the Fourier language makes them increasingly systematic, reusable and general. In the terminology of Section~\ref{sec:framework}, the relevant assertion is that, on an appropriate common comparison scale, some of the corresponding directed transfer costs decrease.

The conceptual distinction is therefore not

\[
\text{no modes}
\longrightarrow
\text{modes},
\]

nor even

\[
\text{no coefficient formula}
\longrightarrow
\text{coefficient formula}.
\]

It is closer to

\[
\text{particular modal and trigonometric procedures}
\longrightarrow
\text{a reusable coordinate architecture}.
\]

Modes cease to function only as distinguished solutions of particular equations and become coordinates through which general data and operators can be organized.

\subsection{Common coordinates and diagonalization}

The significance of the new architecture becomes clearer when the same spatial coordinates are used for different equations. Consider the heat equation

\[
u_t=\kappa u_{xx},
\qquad
u(0,t)=u(\ell,t)=0,
\qquad
u(x,0)=f(x).
\]

If

\[
f(x)
=
\sum_{n=1}^{\infty}
b_n
\sin\left(\frac{n\pi x}{\ell}\right),
\]

then, for suitable data,

\[
u(x,t)
=
\sum_{n=1}^{\infty}
b_n
e^{-\kappa(n\pi/\ell)^2t}
\sin\left(\frac{n\pi x}{\ell}\right).
\]

For the vibrating string

\[
u_{tt}=c^2u_{xx},
\]

with zero initial velocity, the same spatial decomposition gives

\[
u(x,t)
=
\sum_{n=1}^{\infty}
b_n
\cos\left(\frac{n\pi ct}{\ell}\right)
\sin\left(\frac{n\pi x}{\ell}\right).
\]

The temporal dynamics are different, but the decomposition of the initial datum is shared. Heat and wave propagation remain distinct problems while becoming accessible through the same representational machinery.

The underlying mechanism is diagonalization:

\[
-\partial_x^2
\sin\left(\frac{n\pi x}{\ell}\right)
=
\left(\frac{n\pi}{\ell}\right)^2
\sin\left(\frac{n\pi x}{\ell}\right).
\]

Consequently, the PDEs reduce to scalar equations for modal coefficients,

\[
\dot u_n
=
-\kappa
\left(\frac{n\pi}{\ell}\right)^2u_n,
\qquad
\ddot u_n
=
-c^2
\left(\frac{n\pi}{\ell}\right)^2u_n.
\]

Here the productive effect of the representation is especially clear. The same spatial decomposition can be reused across different evolution equations, and the differential operator becomes multiplication by a spectral parameter. This is stronger evidence of rewiring than the mere existence of a coefficient formula: one representation now mediates transfers among distinct mathematical problems.

\subsection{The modern Fourier architecture}

The later development of Fourier analysis makes the scale and persistence of this reorganization still clearer. For \(f,g\in\mathcal S(\mathbb R)\), let

\[
\widehat f(\xi)
=
\int_{\mathbb R}
e^{-ix\xi}f(x)\,dx.
\]

Then

\[
\widehat{\partial_x f}(\xi)
=
i\xi\widehat f(\xi),
\qquad
\widehat{(-\Delta)f}(\xi)
=
|\xi|^2\widehat f(\xi),
\]

while

\[
\widehat{f*g}
=
\widehat f\,\widehat g,
\qquad
\widehat{f(\,\cdot-a\,)}(\xi)
=
e^{-ia\xi}\widehat f(\xi).
\]

These are exact transformations of operations, not merely analogies. Problems involving differential operators become systematically connected with multiplication operators,

\[
P_{\mathrm{diff}}
\longleftrightarrow
P_{\mathrm{multop}},
\]

while convolution problems become connected with pointwise products,

\[
P_{\mathrm{conv}}
\longleftrightarrow
P_{\mathrm{prod}}.
\]

The importance of these correspondences for the present argument is precisely that they are not confined to the vibrating string or the heat equation. The Fourier language becomes reusable across harmonic analysis, spectral theory, partial differential equations and many related areas. A representational architecture originating in a comparatively specific family of problems acquires a much larger transfer domain.

\subsection{Testing the Rewiring Thesis}

The finite family introduced in Section~\ref{sec:string},

\[
F_{\mathrm{string}}
=
\left\{
P_{\mathrm{prop}},
P_{\mathrm{mode}},
P_{\mathrm{repr}},
P_{\mathrm{coeff}},
P_{\mathrm{lim}}
\right\},
\]

provides a useful local test of the formalism, but it should not be treated as an empirical measurement of Fourier's historical breadth.

The strongest claim we need is that the Fourier architecture substantially reorganizes the relation among

\[
P_{\mathrm{repr}},
\qquad
P_{\mathrm{coeff}},
\qquad
P_{\mathrm{mode}}.
\]

Suppose that an independently chosen comparison scale and significance threshold \(\varepsilon\) classify the changes in

\[
(P_{\mathrm{repr}},P_{\mathrm{coeff}})
\qquad\text{and}\qquad
(P_{\mathrm{coeff}},P_{\mathrm{mode}})
\]

as substantial. Then

\[
B_{\mathrm{string}}
=
\left\{
(P_{\mathrm{repr}},P_{\mathrm{coeff}}),
(P_{\mathrm{coeff}},P_{\mathrm{mode}})
\right\}
\subset
A_{\mathcal I_{\mathcal F}}
(F_{\mathrm{string}},\varepsilon),
\]

and

\[
V(B_{\mathrm{string}})
=
\left\{
P_{\mathrm{repr}},
P_{\mathrm{coeff}},
P_{\mathrm{mode}}
\right\}.
\]

Hence

\[
\frac{|V(B_{\mathrm{string}})|}
{|F_{\mathrm{string}}|}
=
\frac35.
\]

This fraction is only a conditional illustration of the breadth criterion.
It depends on the reference family, comparison scale and
significance threshold, and its numerical value would change if the problem
family were refined or coarsened. Hence it should  not be interpreted as an
intrinsic quantitative measure of Fourier's historical significance. The
purpose of \(F_{\mathrm{string}}\) is more modest: to show that the formalism
can register a non-local deformation across several independently specified
problems rather than a mere change in the status of a single problem.

The evidence of breadth is broader than this finite witness.
The same trigonometric coordinates work both for heat and wave evolution;
differential operators become multiplier problems; convolution becomes
multiplication; and the architecture propagates into problem families that
were not part of the original string controversy. These transfers provide the
stronger historical basis to regard Fourier analysis as a large-scale
reorganization of mathematical practice.

Persistence should be understood qualitatively in the same historical sense.
Once Fourier coordinates, spectral decompositions and the associated operator
correspondences were assimilated, they did not remain local devices tied to the
original problems but became reusable components of analysis. In the
terminology of Section~\ref{sec:framework}, this provides evidence for a
persistent core.

The innovation is also expansive. Once frequency becomes an independent
mathematical variable, new questions become natural: how regularity is encoded
by spectral decay, where a function is localized in frequency, which functions
\(m(\xi)\) define useful Fourier multipliers, and how strongly a function can
be localized simultaneously in physical and frequency variables. These lead
naturally to regularity theory, multiplier theory, uncertainty principles,
frequency decompositions and time-frequency analysis. They are not merely the
original string or heat questions rewritten in another notation. In the
terminology of Section~\ref{sec:framework},

\[
\mathcal N_{\mathcal I_{\mathcal F}}
\neq
\varnothing.
\]

\subsection{What the representation reveals and hides}

The triangular profile gives a concrete example of what the new representation reveals. In the symmetric case, Section~\ref{sec:string} gave

\[
|b_n|
=
\frac{8h}{\pi^2n^2},
\qquad
n\ \text{odd}.
\]

In physical space, the corner appears as a failure of differentiability; in spectral variables, the same irregularity is encoded in the asymptotic decay of the coefficients. The singularity has not disappeared: instead, its mathematical form has changed.

But Fourier representation is not a universal simplification. Spatial localization, immediate in physical variables, becomes distributed across frequencies. Likewise, although convolution becomes multiplication,

\[
\widehat{f*g}
=
\widehat f\,\widehat g,
\]

pointwise multiplication becomes convolution,

\[
\widehat{fg}(\xi)
=
\frac{1}{2\pi}
(\widehat f*\widehat g)(\xi)
\]

with the normalization we chose.

Operations that become simple in one direction may become non-local in another. This is exactly why the framework uses deformation rather than contraction: a productive representation redistributes difficulty.

The general theory of such affordances and blind spots will be developed in Section~\ref{sec:blindspots}. For the present purpose, the conclusion is enough. Fourier's importance is not captured merely by the existence of trigonometric series, coefficient formulas or particular PDE solutions. The deeper transformation is the stabilization and propagation of an operational architecture that reorganizes transfers across many mathematical problems.

Within the qualitative terminology fixed in Section~\ref{sec:framework}, Fourier analysis remains our paradigmatic historical case of conceptual rewiring. The finite witness above illustrates how the formalism records part of that change; it is not intended as a numerical proof of its magnitude.

Before Fourier, many of the relevant objects, methods and solutions already existed. After the Fourier language became operational, they could increasingly be treated as parts of a common mathematical architecture.

This is the sense in which Fourier rewired the problem space.

%%%%%%%%%%%%%%%%%%%%%%%%%%%%%%%%%%%%%%%%%%%%%%%
%%%%%%%%%%%%%%%%%%%%%%%%%%%%%%%%%%%%%%%%%%%%%%%
%%%%%%%%%%%%%%%%%%%%%%%%%%%%%%%%%%%%%%%%%%%%%%%

\section{Musical notation as productive representation}
\label{sec:notation}

Fourier analysis provides a particularly clean mathematical example of productive representation. We now move to a very different domain, not because mathematics and music obey the same rules, but because music offers another historical case in which a representation does more than preserve what already exists.

The comparison must remain limited. A musical score is neither a proof nor a complete encoding of a performance, and we will not impose the formalism of Section~2 literally on musical practice. Music serves a different purpose here: to test independently the mechanism described by the Productive Representation Principle. Can a representation change which operations are naturally available, persist beyond the circumstances that produced it, and make new practices possible?

\subsection{From memory to an operational representation}

Music existed long before the musical notation, and sophisticated musical cultures have flourished without anything resembling the modern Western score. Notation is therefore neither a necessary condition for complex music nor a universal endpoint of musical development. Our claim is narrower. Within the Western tradition, notation developed over centuries from systems strongly dependent on oral practice and memory toward increasingly autonomous symbolic representations of musical structure \cite{Grier2021}.

Early neumatic notation could support memory and indicate aspects of melodic movement, but later developments in staff notation made pitch relations increasingly explicit and allowed written music to function with less dependence on prior knowledge of the melody. Schematically, one moves from

\[
\text{notation}+\text{memory of the music}
\longrightarrow
\text{performance}
\]

toward the stronger possibility

\[
\text{notation}
\longrightarrow
\text{reconstruction of previously unknown music}.
\]

The important change is not merely an increase in precision. Relations that previously had to be supplied by memory or oral transmission become accessible through the representation itself.

Rhythmic notation makes the reciprocal relation between representation and practice even clearer. Busse Berger describes how increasingly precise systems of rhythmic notation developed partly in response to the growing complexity of medieval polyphony, but also how those notational resources subsequently enabled new compositional constructions \cite{BusseBerger2002}. The relation is therefore not simply

\[
\text{music}
\longrightarrow
\text{notation},
\]

but rather

\[
\text{musical practice}
\longrightarrow
\text{notation}
\longrightarrow
\text{new musical practice}.
\]

Once stabilized, the representation becomes part of the machinery through which music itself can be organized.

\subsection{Making time spatial}

One of the deepest effects of notation is also one of the most familiar: it turns aspects of musical time into spatial relations,

\[
\text{temporal musical process}
\longrightarrow
\text{spatial symbolic object}.
\]

A performance unfolds sequentially. Written notation does not abolish memory, embodied practice or convention, but it allows events separated in performance time to be simultaneously available to be inspected. In a polyphonic score, simultaneities are aligned vertically while temporal development is displayed horizontally.

This changes the operations that can naturally be performed on musical material. Distant passages can be compared without being reperformed; several voices can be aligned; patterns can be copied, displaced or transposed; and one component can be modified while another is held fixed. The score becomes, in this sense, a \emph{workspace}.

This is close to the proposal of Colyvan and Hall that mathematical and musical notations may function as models rather than merely as neutral records \cite{ColyvanHall2026}. Our claim is more specific. Notation becomes productive when the represented object supports transformations, comparisons and combinations that are substantially harder to control directly in the temporal process itself. In the language of Principle~\ref{pr}, representation changes the available operations.

The representation does not create melody, rhythm or polyphony from nothing. Also, it does not replace the musical practices from which it arose. Its effect is to reorganize the environment in which those structures can be inspected and manipulated. Representation therefore participates in musical production rather than merely following it.

\subsection{Productive without being complete}

The musical case also makes clear that productivity does not require informational completeness. A score stabilizes pitch relations, metric organization, nominal durations, simultaneities and large-scale structural relations, while leaving many dimensions of performance only partially specified. Timing, articulation, timbre, balance, phrasing and other aspects of realization depend on conventions and practices that are not exhausted by the written object.

This incompleteness is not simply a defect. Selectivity is part of what makes notation operational. By suppressing many details of individual realization, the score makes other relations stable enough to be copied, compared, transformed and transmitted. A representation can therefore become more useful for certain operations precisely because it does not preserve every distinction of its object.

This point also guards against treating Western staff notation as a neutral or universal description of music. Notational categories emerge within particular musical practices and may not organize other traditions equally well. Schuiling has accordingly argued that notation should be understood as a cultural practice rather than as a transparent encoding of an independently given musical object \cite{Schuiling2019,Schuiling2022}.

For the present argument, the conclusion is simple:

\[
\text{productive representation}
\not\Rightarrow
\text{complete representation}.
\]

A representation can be conceptually generative because of the distinctions it stabilizes and the operations those distinctions support. The complementary question about which distinctions are lost, collapsed or badly organized in the process will be treated formally in the next section.

\subsection{A distributed transformation}

Musical notation also illustrates that large conceptual transformations need not have a unique inventor or a sharply defined date. The Western notational system developed across centuries through the work of theorists, scribes, teachers, composers, performers and many anonymous practitioners. Guido d'Arezzo occupies an important place in the history of pitch notation and pedagogy, but neither staff notation nor rhythmic notation emerged in a single conceptual act \cite{BusseBerger2002, Grier2021}.

A schematic trajectory is

\[
\begin{aligned}
\text{oral practice}
&\longrightarrow
\text{partial mnemonic representation}\\
&\longrightarrow
\text{more explicit pitch and rhythmic representation}\\
&\longrightarrow
\text{spatial coordination of polyphony}.
\end{aligned}
\]

This should not be read as a universal or strictly linear history of music. It describes one trajectory within the Western notational tradition. Its relevance is structural: relations previously dependent on memory, convention or real-time coordination become progressively available as explicit symbolic relations, and the cumulative effect exceeds improved preservation.

Musical notation therefore exhibits, outside mathematics, the mechanism we have called productive representation. Its effects are substantial, distributed across many musical activities, historically persistent and capable of enabling practices that were difficult to control before the relevant symbolic resources existed. We do not claim that this makes musical notation a conceptual revolution in the formal mathematical sense of Section~2. Music is not evidence for the Rewiring Thesis. It is an independent test of the mechanism underlying the Productive Representation Principle.

The comparison with Fourier can consequently be stated without forcing an analogy. In Fourier analysis,

\[
\text{representation}
\longrightarrow
\text{new transfers among mathematical problems};
\]

in musical notation,

\[
\text{representation}
\longrightarrow
\text{new operations on musical structure}
\longrightarrow
\text{new musical practices}.
\]

In both cases, a representation becomes more than a description of what already exists. It becomes an environment in which previously difficult operations become ordinary and new possibilities become natural.

Music, however, makes the price of that productivity especially visible. The same act of representation that stabilizes some distinctions necessarily leaves others less accessible. A score is not a performance, and that difference is not mere noise. We therefore turn from the productivity of representations to the complementary question: what, exactly, do they leave out?

%%%%%%%%%%%%%%%%%%%%%%%%%%%%%%%%%%%%%%%%%%%%%%%%%
%%%%%%%%%%%%%%%%%%%%%%%%%%%%%%%%%%%%%%%%%%%%%%%%%
%%%%%%%%%%%%%%%%%%%%%%%%%%%%%%%%%%%%%%%%%%%%%%%%%

\section{What representations leave out}
\label{sec:blindspots}

The previous section showed that a representation can be productive without being complete. We now separate two different questions that should not be conflated. The first concerns \emph{information}: which distinctions survive a representation? The second concerns \emph{organization}: which operations and tasks become natural in the represented variables?

The distinction is important. Fiber structure can characterize information loss, but it cannot by itself explain operational productivity. Two injective representations have the same trivial fiber structure and may nevertheless support very different mathematical operations. Conversely, a non-injective representation may be highly productive for tasks that do not require the distinctions it suppresses.

Thus representational affordances and representational blind spots are related, but they are not reducible to a single formal property.

\subsection{Fibers and failure of exact recovery}

Let

\[
R:X\longrightarrow Z
\]

be a representation. For \(x\in X\), its fiber is

\[
R^{-1}(R(x))
=
\{y\in X:R(y)=R(x)\}.
\]

Elements in the same fiber are indistinguishable from the point of view of \(R\).

Now let

\[
T:X\longrightarrow Y
\]

be a quantity, property or task whose exact value is to be recovered from the representation. Then exact recovery is possible precisely when \(T\) is constant on the fibers of \(R\).

Equivalently, there exists a map

\[
\widetilde T:R(X)\longrightarrow Y
\]

such that

\[
T=\widetilde T\circ R
\]

if and only if

\[
R(x_1)=R(x_2)
\quad\Longrightarrow\quad
T(x_1)=T(x_2).
\]

Consequently, if

\[
R(x_1)=R(x_2),
\qquad
T(x_1)\neq T(x_2),
\]

then \(T\) cannot be recovered exactly from \(R(x)\) alone.

This is a precise sense in which a limitation can be representational rather than computational. No increase in computation performed solely on \(R(x)\) can reconstruct a distinction that \(R\) has collapsed.

The qualification ``exactly'' matters. A representation may still support probabilistic inference, approximation or statistically useful prediction of \(T\). The statement concerns exact recovery for all inputs.

The criterion is also task-dependent. Two elements of \(X\) may be legitimately identified for one family of tasks and crucially different for another. Abstraction itself depends on suppressing distinctions; the relevant question is whether the suppressed distinctions matter for the task under consideration.

\subsection{Operations and productive affordances}

The previous subsection concerns information loss. Operational productivity is a different issue.

Suppose

\[
A:X\longrightarrow X
\]

is an operation and there exists

\[
B:Z\longrightarrow Z
\]

such that

\[
R\circ A
=
B\circ R.
\]

Then the represented effect of \(A\) can be obtained entirely inside the represented space:

\[
R(x)
\longmapsto
B(R(x))
=
R(Ax).
\]

As emphasized in Section~\ref{sec:framework}, this intertwining relation is not by itself sufficient for productivity. A constant representation can satisfy such a relation trivially. The representation becomes operationally productive only when the induced operation \(B\) preserves structure relevant to the task and makes useful operations, comparisons or transfers substantially easier.

Thus two distinct questions arise:

\[
\begin{aligned}
\text{Information question:}\quad
&\text{which distinctions does }R\text{ preserve?}\\[3pt]
\text{Operational question:}\quad
&\text{how does }R\text{ reorganize useful operations?}
\end{aligned}
\]

The first is controlled by fibers. The second depends on the geometry, algebra or organization of the represented space.

This distinction is essential for Fourier analysis. The Fourier transform, on suitable function spaces, is invertible and therefore does not derive its productivity from collapsing information. Its power comes from reorganizing operations: differentiation becomes multiplication, convolution becomes pointwise product, and spectral localization becomes explicit.

A productive affordance should therefore be understood as an operational advantage created by the organization of a representation, not merely as a consequence of its information content.

\subsection{The score and the performance}

Music makes the information side of the distinction especially concrete. Let \(P\) denote a performance and let

\[
R(P)=S
\]

be its written score. Even in an idealized setting, the map from performances to scores is highly non-injective. Distinct performances may satisfy

\[
R(P_1)=R(P_2)=S
\]

while differing in timing, dynamics, articulation, balance, pedalling, timbre or phrase shaping. Much of the work of Widmer and collaborators on computational performance analysis concerns precisely such dimensions of played music \cite{Widmer2003,WidmerGoebl2004}.

The score therefore places many musically distinct performances in the same fiber. Yet this does not make it a poor representation. Its selectivity is one source of its usefulness. By abstracting from many details of individual realization, notation stabilizes pitch relations, metric organization, simultaneity and large-scale structure, making those features available for comparison, transformation and transmission.

A recent computational example illustrates the same point. Peter, Hu and Widmer study real-time score following by converting an incoming piano performance into a symbolic note representation and aligning that representation with the score \cite{PeterHuWidmer2025}. The intermediate representation discards acoustic information, yet it may be better adapted to the alignment task. More retained information is therefore not automatically equivalent to a more useful representation.

The limitation appears when the task changes. The question

\[
\text{Where are we in the score?}
\]

depends largely on distinctions preserved by the symbolic representation, whereas

\[
\text{How is this passage being performed?}
\]

depends on distinctions that the same representation suppresses.
There is therefore no task-independent notion of an optimal representation.

Related difficulties arise under changes of domain. Marták, Hu and Widmer show that automatic music-transcription systems may deteriorate under shifts in sound and musical distribution, including changes in genre, dynamics and polyphonic complexity \cite{MartakHuWidmer2026}. Such behavior reminds us that representations and the systems built around them encode assumptions about which distinctions matter and under which conditions those distinctions remain useful.

\subsection{Loss and misorganization}

Non-injectivity describes one important failure mechanism, but not the only one.

The first is \emph{loss}. A representation identifies objects that a relevant task must distinguish:

\[
R(x_1)=R(x_2),
\qquad
T(x_1)\neq T(x_2).
\]

The distinction required by \(T\) is then unavailable from \(R(x)\) alone.

The second is \emph{misorganization}. A representation may preserve the relevant information while arranging it in variables or categories that fail to expose the structure needed for a task. In this case the obstruction is not informational loss. The difficulty lies in how the available information has been organized.

Thus

\[
\boxed{\text{loss}\neq\text{misorganization}}.
\]

The distinction matters because the natural responses differ. Loss may require enriching a representation. Misorganization may require reorganizing or replacing it even when no information has been discarded.

Musical representation again provides a concrete example. Western staff notation naturally privileges categories such as discrete notes, pitches, durations, measures and vertically coordinated events. These categories are highly productive within the traditions in which the notation developed, but they need not organize every musical practice equally well.

Nuttall, Serra and Pearson study the different forms taken by a \emph{svara} in Carnatic performance and emphasize the importance of melodic context and ornamentation in its realization \cite{NuttallSerraPearson2025}. In such a setting, representing the relevant musical event primarily as a context-independent discrete pitch may obscure precisely the structure one wishes to understand. The problem is not necessarily that too little information is retained; the chosen variables may organize the phenomenon poorly.

Schematically, suppose

\[
R(M)=(z_1,\ldots,z_N)
\]

represents a rich process \(M\). Even if the map \(R\) preserves extensive information, the relevant mathematical or musical structure may not align naturally with the coordinates \(z_j\). Estimating those coordinates with increasing accuracy need not solve the conceptual problem. One may simply be manipulating poorly chosen variables more precisely.

The same issue persists for learned representations. Learning coordinates from data does not make representation neutral: the learned variables still organize some regularities more effectively than others. Ibáñez-Martínez, Nkama, Poltronieri, Serra and Rocamora illustrate this in controllable music generation, where intended semantic decompositions of learned latent spaces need not coincide with the information those components actually contain \cite{IbanezMartinezEtAl2026}. Hence

\[
\text{representation label}
\neq
\text{represented structure}.
\]

A representation must therefore be understood both through the distinctions it preserves and through the organization it imposes on those distinctions.

\subsection{Affordances, blind spots, and representational revision}

We can now sharpen the terminology.

A representation has an \emph{affordance} relative to a family of tasks when its organization makes a relevant operation, comparison or transfer substantially easier.

A representation has an \emph{informational blind spot} relative to a task when it collapses distinctions required for exact recovery of that task.

It may also be \emph{misorganized} relative to a family of tasks when the relevant information is present but arranged in variables or categories that make the needed structure difficult to expose.

These phenomena not necessarily coincide. An invertible representation can have powerful affordances and severe organizational disadvantages without losing information. A non-injective representation can be highly effective whenever the distinctions it suppresses are irrelevant to the target tasks.

This also clarifies the relation with transfer structure. There is no requirement that

\[
c_{L_{\mathcal I}}
\bigl(
\iota_{\mathcal I}(P),
\iota_{\mathcal I}(Q)
\bigr)
<
c_L(P,Q)
\]

for every pair \(P,Q\). A productive representation may lower some transfer costs while increasing others. Conceptual rewiring concerns the resulting pattern of deformation, not universal improvement.

A more serious situation arises when several important tasks repeatedly expose a common representational obstruction. In the case of information loss, one may repeatedly encounter

\[
R(x_1)=R(x_2),
\qquad
T_j(x_1)\neq T_j(x_2).
\]

In the case of misorganization, several problems may remain difficult because the variables supplied by \(R\) systematically fail to expose a common structure.

At that point, greater computational sophistication inside the existing represented space may not address the common source of difficulty. The representation itself becomes a candidate for revision.

One may then seek

\[
R':X\longrightarrow Z'
\]

whose variables preserve or expose the distinctions relevant to the problematic tasks and whose organization supports a more effective family of operations.

The significant upgrade is not merely to replace one encoding by another: it is instead to diagnose a recurrent obstruction as specifically representational and to construct a new representation that reorganizes several related tasks.

This is the motivation for the \emph{Representational Escape Problem} developed in Section~\ref{sec:ai}.

Representational revision, however, is only one possible route to conceptual rewiring. A mathematical innovation may instead introduce a new structural mediator, objectify something previously treated as background, enlarge the admissible universe, or combine several such mechanisms. These categories are not necessarily mutually exclusive.

The next section therefore asks whether the Rewiring Thesis survives cases whose dominant mechanism is not most naturally described as representational change.

%%%%%%%%%%%%%%%%%%%%%%%%%%%%%%%%%%%%%%%%%%%%%%%%%
%%%%%%%%%%%%%%%%%%%%%%%%%%%%%%%%%%%%%%%%%%%%%%%%%
%%%%%%%%%%%%%%%%%%%%%%%%%%%%%%%%%%%%%%%%%%%%%%%%%

\section{Stress tests for the Rewiring Thesis}
\label{sec:stress-tests}

Fourier analysis is an unusually favorable example for the framework developed in this paper. Its representation is explicit, the relevant operations can be written by exact formulas, and the persistence of the resulting transfers is visible throughout modern analysis. Precisely for that reason, it leaves an important question open: is the Rewiring Thesis genuinely broader than successful representational change?

If so, comparable reorganizations should also occur when the decisive innovation is not naturally described as the replacement of one representation by another. We therefore consider four deliberately different cases. The aim is not to reconstruct their complete histories or to assign numerical transfer costs. We ask instead which mechanism reorganizes the relevant problem space and whether the effect extends beyond the problem that initially motivated it.

\subsection{Galois: from equations to structural mediation}

The classical problem of solving polynomial equations provides a particularly clean first test. After the success of explicit formulas for equations of degrees three and four, the corresponding question for higher degrees was naturally posed in constructive terms: can the roots of a polynomial be expressed from its coefficients using an allowed collection of algebraic operations?

Galois changed the organization of this problem. His work, written in the early 1830s and published by Liouville in 1846, related solvability by radicals to the permutations of the roots that preserve their algebraic relations \cite{Galois1846}. In modern language, if $f$ is a polynomial over a field $K$ of characteristic zero and $K_f$ is a splitting field, one associates the Galois group

\[
G_f=\operatorname{Gal}(K_f/K),
\]

and the classical solvability question becomes the structural criterion

\[
\text{$f$ is solvable by radicals}
\quad\Longleftrightarrow\quad
\text{$G_f$ is a solvable group}.
\]

This modern formulation uses concepts and terminology that were stabilized only after Galois, so it should not be projected back unchanged onto his original work. That later stabilization is nevertheless part of what matters here. Earlier work  by Lagrange and Abel had already moved beyond a simple search for increasingly ingenious formulas; Galois introduced a particularly powerful new structural mediator through which solvability questions for different equations could be organized and compared.

Schematically,

\[
\text{polynomial equations}
\longrightarrow
\text{symmetries of roots}
\longrightarrow
\text{group structure}
\longrightarrow
\text{solvability}.
\]

The Galois group is therefore not only another representation of the roots. It is a \emph{structural mediator}: equations that look very different at the level of coefficients can be compared through the structures of their associated groups. The new object also generates questions that are no longer subordinate to the original problem of finding formulas for roots.

This gives a first mechanism of rewiring distinct from Fourier's. Fourier reorganizes problems through productive coordinates; Galois introduces a new structural object through which previously separate problems become transferable.

\subsection{Riemann: when the background becomes an object}

Riemann illustrates a second mechanism. Building on earlier developments in intrinsic and non-Euclidean geometry, his 1854 habilitation lecture, published posthumously in 1868, substantially generalized the range of geometries that could be treated mathematically \cite{Riemann1868}. The geometry of the ambient space itself became variable mathematical structure.

In modern language, one considers a smooth manifold $M$ equipped locally with a metric tensor

\[
g
=
\sum_{i,j}
g_{ij}(x)\,
dx^i\otimes dx^j.
\]

Lengths, angles, geodesics and curvature are determined by $g$, while the underlying geometry is independent of the particular coordinates used to describe it. The decisive conceptual step is therefore not the introduction of a more efficient coordinate system. Coordinates become auxiliary descriptions of an object whose geometry can itself vary.

The effect on the problem space is substantial. Local differential information becomes connected with curvature and other geometric invariants; differential equations become natural tools for studying geometric structures; and local properties become related to global questions about the underlying space. At the same time, questions arise that make no sense when geometry is treated as a fixed background: which metrics exist on a given manifold, which curvature conditions they satisfy, and how geometric structures behave under deformation.

The mechanism can be described as \emph{objectification}: something that previously functioned as part of the background of a problem becomes a variable mathematical object that can itself be classified, transformed and compared.

This is conceptually different from the Fourier case. Fourier changes the coordinates through which existing objects are manipulated; Riemann changes what counts as an object of investigation. Rewiring can therefore arise from a transformation of mathematical ontology, not only from a new representation.

\subsection{Lebesgue and Schwartz: enlarging the admissible universe}

A third mechanism appears when the obstruction lies in the class of objects admitted by a theory. Instead of finding a better method inside a fixed universe, mathematics may enlarge that universe until operations that were previously exceptional become stable.

Lebesgue's theory of measure and integration provides a fundamental example. His 1902 PhD work introduced an integration theory adapted to a substantially broader class of functions and, crucially, to limiting procedures that are difficult to control within the Riemann framework \cite{Lebesgue1902}. Its conceptual effect is therefore larger than the ability to integrate additional functions. Measure-theoretic language provides common environments in which convergence, approximation, Fourier analysis, probability and differential equations can interact.

The later development of spaces such as $L^p$ makes the reorganization especially visible. Pointwise behavior ceases to be the only natural way to compare functions; norm convergence, almost-everywhere equivalence and integrability classes become mathematical structures in their own right. New transfers become possible because questions that previously required different languages can now be posed inside a common functional framework.

Schwartz's theory of distributions gives an even sharper example of the same mechanism. Let

\[
\mathcal D(\Omega)=C_c^\infty(\Omega),
\]

and let

\[
\mathcal D'(\Omega)
=
\mathcal D(\Omega)'
\]

denote its continuous dual. For $T\in\mathcal D'(\Omega)$, differentiation is defined by

\[
\langle \partial_j T,\varphi\rangle
=
-
\langle T,\partial_j\varphi\rangle.
\]

Every distribution therefore has derivatives of all orders in the distributional sense, independently of whether it is represented by a classically differentiable function. Objects such as the Dirac mass enter the same analytical language as ordinary functions, and differential equations without classical solutions may acquire weak or distributional solutions \cite{Schwartz1950}.

The structural mechanism is

\[
\begin{aligned}
\text{classical objects}
&\longrightarrow
\text{enlarged function spaces}\\
&\longrightarrow
\text{extended operations}
\longrightarrow
\text{new notions of solution}.
\end{aligned}
\]

Here expansion is not merely a side effect of the innovation: it is the mechanism itself. By changing the admissible universe, operations and solutions that were previously unavailable become ordinary parts of the theory.

Lebesgue and Schwartz therefore provide a third form of rewiring: not a new coordinate system and not primarily a new mediator, but an enlargement of the class of objects on which mathematical operations are allowed to act.

\subsection{Wiles: separating status from architecture}

The previous examples are positive cases of structural reorganization. Fermat's Last Theorem serves a different purpose: it helps separate the solution of a major problem from the mathematical architecture through which that solution becomes possible.

The theorem asserts that

\[
x^n+y^n=z^n,
\qquad
n\in\mathbb N,\quad n>2,
\]

has no nonzero integer solutions. Wiles's 1995 paper proved the modularity of semistable elliptic curves over $\mathbb Q$, from which Fermat's Last Theorem follows through the previously established connection involving the Frey curve and Ribet's theorem \cite{Ribet1990, Wiles1995}. Together with the companion paper of Taylor and Wiles, the argument developed deformation-theoretic methods relating Galois representations and Hecke algebras \cite{TaylorWiles1995}.

From the viewpoint of the status map introduced in Section~2, the most visible outcome is simply

\[
\sigma_L(\mathrm{FLT}):
\text{open}
\longrightarrow
\text{solved}.
\]

But the route to that status change runs through an architecture remote from the original elementary-looking Diophantine equation:

\[
\text{Diophantine equation}
\longrightarrow
\text{elliptic curve}
\longrightarrow
\text{Galois representations and modular forms}.
\]

It would be historically and conceptually misleading to attribute this entire architecture to Wiles. The Frey curve, the modularity conjecture, Serre's ideas and Ribet's theorem had already created crucial transfers between Diophantine equations, elliptic curves and modular forms. Wiles and Taylor--Wiles completed a decisive part of this architecture and introduced methods whose significance extends beyond Fermat's equation.

The example therefore separates two levels that a theory of conceptual change should not conflate. The statement

\[
\text{open}\longrightarrow\text{solved}
\]

records a change in the status of one problem. The network connecting elliptic curves, modular forms, Galois representations and deformation theory records changes in transfer structure. Both are present in the historical episode, but they are not the same achievement.

This is why Fermat's Last Theorem is useful here. The Rewiring Thesis does not claim that structural change is more valuable than solving a monumental problem, nor that every great theorem must reconfigure a field. It distinguishes the kinds of change involved. A proof can be extraordinary because of the problem it settles, while the architecture surrounding the proof can independently alter relations across a much broader domain.

\subsection{What survives the tests}

Taken together with Fourier, the four stress tests reveal several distinct mechanisms of conceptual reorganization:

\[
\begin{array}{rcl}
\text{Fourier} &:& \text{productive representation},\\
\text{Galois} &:& \text{structural mediation},\\
\text{Riemann} &:& \text{objectification of background},\\
\text{Lebesgue--Schwartz} &:& \text{enlargement of admissible objects}.
\end{array}
\]

The Wiles case plays a complementary role: it shows why status change and structural reorganization must be distinguished even when they occur within the same mathematical episode.

The common feature of the positive cases is not a particular technique or representation. It is a persistent change in the relations among problems. An innovation outlives the question that initially motivated it, creates reusable routes between previously separated problems and, in the strongest cases, alters which objects and questions are naturally available.

The stress tests therefore support the logical hierarchy of the paper. The Productive Representation Principle identifies one powerful mechanism for conceptual change, but the Rewiring Thesis is broader. Mathematical problem spaces may also be reorganized by new mediating structures, by turning background into object, by enlarging the admissible universe, or by mechanisms not represented in these historical examples.

This distinction matters when we turn to artificial intelligence. Representational escape will provide one concrete test of conceptual creativity, but it cannot exhaust the question. The broader issue is whether an artificial system could autonomously produce any substantial, non-local and persistent reorganization of mathematical transfer structure.

We now turn to that question.

%%%%%%%%%%%%%%%%%%%%%%%%%%%%%%%%%%%%%%%%%%%%%%%%%
%%%%%%%%%%%%%%%%%%%%%%%%%%%%%%%%%%%%%%%%%%%%%%%%%
%%%%%%%%%%%%%%%%%%%%%%%%%%%%%%%%%%%%%%%%%%%%%%%%%

\section{Artificial intelligence and representational escape}
\label{sec:ai}

The distinctions developed in the previous sections become especially relevant for artificial intelligence. The recent discussions of AI and mathematics are naturally organized around capability: which problems a system can solve, how difficult they are, how much human intervention is required, whether the resulting proofs are correct, and whether mathematicians can understand them. These questions are fundamental, but {\bf they do not exhaust mathematical creativity}.

Tao has suggested that, rather than waiting for agreement about the eventual capabilities of artificial systems, one may condition on the possibility that they become able to perform a significant fraction of research-level mathematical work and ask which mathematical values would then remain to be optimized \cite{Tao2026}. Commelin, Jamnik, Ochigame, Taelman and Venkatesh likewise emphasize the need to examine how AI may affect mathematical practice, values and infrastructure \cite{CommelinEtAl2026}. We adopt the same conditional attitude. Suppose that an artificial system can solve difficult open problems, verify its proofs and communicate them successfully. What additional evidence would indicate conceptual creativity in the sense developed here?

The Rewiring Thesis provides one answer. Exceptional problem solving may occur entirely within objects, representations and transfer structures inherited from existing mathematics. Conceptual revolution requires something different: a substantial, non-local and persistent reorganization of relations among problems.

This focus is compatible with, but distinct from, other recent discussions of mathematical innovation and AI. Jutla and Sharma also connect the vibrating-string controversy and Fourier analysis with artificial mathematical creativity, emphasizing the role of patterns supplied by nature \cite{JutlaSharma2026}. Our question is different. We ask whether an artificial system can diagnose that a representation inherited from existing mathematics has itself become part of the obstruction and autonomously replace it by another one.

\subsection{Representation learning, invention, and revolution}

Modern artificial systems already learn representations. Schematically,

\[
R_\theta:X\longrightarrow Z,
\]

where the parameters \(\theta\) are adjusted from data or feedback so that the latent space \(Z\) supports one or more tasks. The resulting coordinates need not resemble variables chosen by a human designer and may encode regularities that were never explicitly supplied. It would therefore be mistaken to contrast humans who invent representations with machines that merely manipulate fixed ones.

Representation learning, however, is not by itself conceptual innovation. A learned encoding may improve performance dramatically while remaining internal to one system, optimized for a particular objective and unavailable as a reusable mathematical language.

It is useful to distinguish three levels:

\[
\text{representation learning}
\longrightarrow
\text{representation invention}
\longrightarrow
\text{representational revolution}.
\]

The arrows indicate increasing conceptual demands, not necessary stages through which every innovation must pass.

By \emph{representation learning} we mean the construction or adaptation of an encoding that improves performance on one or more tasks.

By \emph{representation invention} we mean the production of an externalizable system of variables, coordinates or representational relations that was not already available as that representation in the relevant body of mathematics. This is a stronger novelty claim than merely introducing a familiar representation into a prompt in which it was not explicitly mentioned.

Several novelty levels should therefore be distinguished:

\[
\begin{aligned}
&\text{new to the prompt},\\
&\text{new to the current system interaction},\\
&\text{new to the target problem family or domain},\\
&\text{new to mathematics}.
\end{aligned}
\]

Only the last supports an unqualified claim of representational invention. The weaker levels may still be mathematically important. In particular, autonomously importing an existing representation into a domain where it was not previously used may produce substantial conceptual reorganization without constituting the invention of that representation.

A \emph{representational revolution} is therefore defined by its consequences rather than by novelty alone. It occurs when the introduction, invention or novel deployment of a representation produces the kind of broad and persistent reorganization described by the Rewiring Thesis.

Thus

\[
\text{representation invention}
\not\Rightarrow
\text{representational revolution},
\]

and conversely a representational revolution need not require that every component of the new language be unprecedented in mathematics.

\subsection{The Representational Escape Problem}

We now return to Question~\ref{q:escape}. Suppose an artificial system works on a family

\[
F=\{P_1,\ldots,P_N\}\subset\mathcal P_L
\]

expressed in a conceptual language \(L\). Several problems remain resistant despite powerful search, deduction, experimentation and access to existing mathematical knowledge.

One possible response is to search more effectively within the existing representation. Another is to infer that several difficulties share a specifically representational source: distinctions needed by the tasks are being collapsed, or the available variables organize the relevant structure poorly.

Section~\ref{sec:blindspots} identified these two possibilities as loss and misorganization.

In the case of loss, one may repeatedly encounter

\[
R(x_1)=R(x_2),
\qquad
T_j(x_1)\neq T_j(x_2),
\]

showing that distinctions required by several tasks are unavailable from the representation alone.

In the case of misorganization, the relevant information may be present, yet the variables supplied by

\[
R:X\longrightarrow Z
\]

fail to expose a common structure needed across several problems.

One response is to construct increasingly sophisticated procedures that compensate for these limitations while continuing to work inside \(Z\). A more radical response is to diagnose the common obstruction as representational and construct

\[
R':X\longrightarrow Z'
\]

whose variables preserve or expose the relevant structure more effectively.

We call this move \emph{representational escape}.

The term is deliberately narrow. Representational escape does not include every case in which a mathematical language changes, a new object is introduced, or the admissible universe is enlarged. It refers specifically to the diagnosis and replacement of an obstructive representation.

Representational escape and representational revolution should also not be identified. Escape names the conceptual move away from a representation judged inadequate for a family of tasks. It becomes revolutionary only if the replacement subsequently produces a sufficiently broad and persistent reorganization of mathematical transfer structure.

\subsection{Prospective evidence for representational escape}

A serious test of representational escape should require more than unusual notation, an opaque embedding or improved benchmark performance. At the moment of discovery, however, one cannot usually establish the full persistence condition of Section~\ref{sec:framework}; persistence is partly historical.

The following criteria should therefore be understood as \emph{prospective evidence} for a beneficial representational rewiring, not as an operational equivalent of the formal definition of conceptual revolution.

Four core criteria are relevant.

\begin{enumerate}

\item \textbf{Representational diagnosis.}
The system identifies a recurrent obstruction shared by several problems and attributes it specifically to the representation in which the relevant objects are being organized: important distinctions are collapsed, or the represented variables systematically fail to expose useful structure.

\item \textbf{Representational replacement.}
The system constructs or selects a substantially different representation

\[
R':X\longrightarrow Z'
\]

without being explicitly told that this is the representation required. The novelty level of the replacement must be stated separately: rediscovery, transfer from another domain and genuine invention are different achievements.

\item \textbf{Beneficial transfer.}
The replacement makes several relevant transfers easier, rather than merely solving the task that motivated its construction. In the language of Section~\ref{sec:framework}, the evidence should include negative transfer deformations for more than one strategically related pair of problems.

\item \textbf{Reuse and prospective persistence.}
The representation remains useful beyond the first successful application, can be externalized and communicated, and supports independent reuse on related problems. Such reuse is evidence toward persistence, although full \(\varepsilon\)-persistence in the sense of Section~\ref{sec:framework} can only be established after subsequent assimilation.

\end{enumerate}

An additional feature is

\begin{enumerate}[resume]

\item \textbf{Expansion.}
The new representation makes natural mathematical objects or questions that were not naturally available in the previous language.

\end{enumerate}

Expansion is not required by the definition of conceptual revolution; it identifies an especially important additional consequence.

These criteria are deliberately stricter in one respect and weaker in another than the general formalism. They are stricter because representational escape is intended to be beneficial: the central transfers should become easier, whereas the general Rewiring Thesis allows substantial deformations of either sign. They are weaker because immediate reuse cannot establish long-term persistence. The criteria therefore identify a prospective subclass of conceptual rewiring rather than reproduce the formal definition.

Nor should they be interpreted as a numerical creativity score. Transfer costs have no universal scale, and the novelty of a representation cannot be inferred from benchmark performance alone. Solving many problems and reorganizing a smaller family through one representation are different achievements.

\subsection{Can representational escape be tested?}

There is an immediate difficulty in turning this idea into an evaluation procedure. Mathematical AI is naturally tested on tasks whose answers are known or whose correctness can be verified. Representational escape has almost the opposite structure. If an evaluator specifies the right representation, the crucial conceptual step has already been supplied. If a benchmark specifies the decisive intermediate object, a system may reach it without diagnosing any inadequacy in the original representation.

One possible experiment is historical reconstruction. A system could be given mathematical knowledge available before a major representational transition, while later terminology and constructions are withheld, and asked to make progress on a carefully selected family of problems.

Such experiments would be informative but would test recovery rather than origination. A model trained on modern mathematical literature may retain traces of the supposedly withheld construction, and reproducing a known historical transition does not establish that the system invented it. At most, it shows that the system can autonomously reconstruct or redeploy the relevant representation under constrained conditions.

A stronger test would therefore use genuinely new problem families. Several tasks would be selected so that substantial progress becomes possible, or dramatically easier, after abandoning a natural initial representation. The system would not be told that representational change is required. It would have to detect a recurrent pattern of failure, propose a replacement and show that the new representation improves transfers across the family.

Evaluation would then occur at two levels. The first is immediate: does the proposed representation produce independently verifiable mathematical gains on several problems? The second is necessarily slower: does the representation remain reusable, generative and capable of entering ordinary mathematical practice?

This second stage is related to what Tao describes as mathematical digestion and canonicalization \cite{Tao2026}. In the terminology of this paper, such assimilation provides evidence for persistence.

The requirement of intelligibility must also be handled carefully. A conceptual innovation need not initially resemble familiar human intuition, and mathematical history contains many structures whose significance became clear only after substantial assimilation. Artificial creativity should therefore not be defined by psychological similarity to human mathematicians.

A purely internal representation raises a different issue. If it cannot be externalized into mathematical variables, relations or operations reusable outside the system, it is unclear in what sense mathematics itself, rather than one computational device, has been reorganized.

The relevant requirement is therefore not immediate transparency but \emph{externalizable structure}. Something produced by the system must be capable of entering shared mathematical practice.

This requirement is relative to the claim being made. A private internal representation may constitute an important computational achievement. Externalizability becomes necessary only when the stronger claim is that the system has reorganized mathematics as a shared body of structures, methods and problems.

This criterion also avoids unnecessary anthropomorphism. Representational escape does not require a claim about whether a machine experiences insight, surprise or understanding. The relevant evidence concerns mathematical consequences. Conversely, a fluent explanation that an existing representation is inadequate is not sufficient if no mathematical organization actually changes.

\subsection{Beyond representation}

Representational escape is only one route to conceptual creativity. The stress tests of Section~\ref{sec:stress-tests} showed that transfer structure can also be reorganized through structural mediation, objectification of what had previously functioned as background, or enlargement of the admissible universe.

These mechanisms should not be understood as mutually exclusive categories. A development can often be described representationally at some level. For example, associating a polynomial with its Galois group defines a map from one class of mathematical objects to another. Nevertheless, describing Galois theory merely as a change of representation would miss what is structurally distinctive in the example: the introduction of a new mathematical object whose internal structure mediates solvability and connects previously separated equations.

These distinctions are therefore meant to identify different explanatory mechanisms, not to place them in separate, mutually exclusive categories. The goal is to understand which mechanism best accounts for the resulting change in transfer structure, while allowing several mechanisms to coexist in the same historical development.

The general AI question is correspondingly broader than representational escape. An artificial system might eventually produce conceptual rewiring by inventing a representation, introducing a new structural mediator, enlarging the admissible class of objects, objectifying a previously fixed background, or through a mechanism for which mathematical history offers no close precedent.

A further possibility, which deserves separate investigation, concerns the
relation between representation and physical intuition. In the
eighteenth-century string problem, mathematical representation largely
followed physical intuition: oscillatory modes gave mathematical form to
patterns that could, at least in principle, be seen or heard. A century after
Fourier, however, Fourier duality became part of the mathematical architecture
of quantum mechanics, where the relevant structures could no longer be
organized by an equally direct classical picture. Representation was no longer
merely encoding an intuition already available to physics; it was helping to
determine the very form in which physical structure could be accessed.

This suggests a stronger version of the representational escape problem.
\textbf{Could an artificial system construct a mathematically productive
representation before a corresponding human intuition exists?} If so, the next
Fourier transform might not merely reorganize problems we already understand.
It might provide the mathematics through which we first learn what there is to
understand.

Representational escape remains especially useful because it isolates one recognizable and testable transition. Artificial systems may become extraordinarily effective navigators of the mathematical world inherited from human mathematics. They may discover routes through structures of remarkable complexity.

But navigation and cartography are different achievements. A navigator finds routes on an existing map; a conceptual innovation may alter the map itself.

Representational escape asks whether a machine can recognize when one of the maps inherited from mathematics has become part of the problem---and construct a better one.

%%%%%%%%%%%%%%%%%%%%%%%%%%%%%%%%%%%%%%%%%%%%%%%%%%
%%%%%%%%%%%%%%%%%%%%%%%%%%%%%%%%%%%%%%%%%%%%%%%%%%
%%%%%%%%%%%%%%%%%%%%%%%%%%%%%%%%%%%%%%%%%%%%%%%%%%
%%%%%%%%%%%%%%%%%%%%%%%%%%%%%%%%%%%%%%%%%%%%%%%%%%

\section{The next Fourier transform}
\label{sec:conclusion}

We began with a vibrating string. By the middle of the eighteenth century, much of the relevant mathematics was already present: propagation formulas, normal modes, increasingly broad notions of admissible functions, and finite-dimensional routes toward the continuum. What was missing was not simply one further calculation. The available descriptions did not yet belong to a common operational language.

Fourier analysis changed that organization. Modes became coordinates, differential operators became multipliers, different evolution equations could be studied through the same spectral decomposition, and questions about regularity, localization and frequency became natural in their own right.

\textbf{The problems were already there. Fourier changed the distances between them.}

This observation motivated our \textbf{Rewiring Thesis}. Some mathematical innovations matter not only because of the problems they solve, but because they reorganize the transfers through which problems can inform one another. The formalism developed in this paper was intended to make that claim
sufficiently precise to distinguish amplitude, breadth and persistence,
while tracking expansion as an additional feature, without pretending
that conceptual change admits a universal numerical metric.

Fourier also revealed one particularly powerful mechanism for such change: productive representation. A representation becomes conceptually productive when it changes the operations available on its objects and thereby alters the relations among problems. Musical notation provided an independent test of this mechanism. Its importance for our argument was precisely that it is both productive and incomplete: its power comes from stabilizing some distinctions while leaving others less explicit. Representation therefore does not become useful by preserving everything. It becomes useful by organizing what can be done.

The mathematical stress tests showed, however, that productive representation does not exhaust conceptual revolution. Transfer structure may also be reorganized by new structural mediators, by turning what had functioned as background into a mathematical object, by enlarging the admissible universe, or by other mechanisms. The Rewiring Thesis is therefore broader than the Productive Representation Principle.

This distinction matters for artificial intelligence. A system may become extraordinarily successful at proving theorems while continuing to operate within the conceptual organization inherited from existing mathematics. That would be genuine mathematical progress, but it would not by itself establish conceptual revolution.

Representational escape identifies one stronger possibility. A system might diagnose that recurrent difficulties arise not merely from insufficient search but from the representation in which a family of problems has been posed, replace that representation, and thereby make previously separate problems accessible through a common structure. Escape, however, is only the conceptual move. It becomes a \emph{representational revolution} if the new structure subsequently produces a broad and persistent reorganization of mathematical transfers. And even representational revolution is only one possible mechanism of conceptual rewiring.

The relevant criterion is therefore structural rather than psychological.
We need not decide whether an artificial system experiences insight or
understanding in a human sense. The mathematical question is what changes
after its contribution. Does a new structure survive the problem that
produced it? Can it be externalized, reused and developed? Do previously
distant problems become transferable? And, in expansive cases, do new
questions become natural?

This is what we mean by the ``next Fourier transform.'' We do not mean another integral transform, and there is no reason to expect the next conceptual revolution to resemble Fourier analysis technically. The phrase denotes an innovation after which a familiar collection of mathematical problems is no longer organized in the same way.

Such transformations need not have a single author or a sharply identifiable moment of birth. The histories considered here repeatedly show that conceptual change may be distributed across people, problems and decades, and that its magnitude may become visible only after the resulting transfers have entered ordinary mathematical practice. For the same reason, the first convincing evidence of artificial conceptual creativity may not be a spectacular isolated theorem.

The more demanding possibility is that an artificial system changes the map itself: distinctions that had structured a subject cease to be fundamental, previously distant problems become manifestations of a common structure, and questions that had not previously been natural become unavoidable.

If that happens, the most significant mathematical contribution of artificial intelligence may not be the first famous open problem it solves. It may be the first time mathematicians realize that, after an idea produced by a machine, problems they had studied for years are no longer where they used to be.

\textbf{That would be the next Fourier transform.}

\subsection*{Acknowledgement}

During the last two years I have had the privilege of collaborating with Gerhard Widmer and Xavier Serra on the organization of the Villa de Bilbao music competition, the Spirio Prize, and related activities connecting mathematics, music and technology. These collaborations have given me the opportunity to discuss music, computation and representation with both of them, often far outside my own field. I have learned a great deal from Gerhard and Xavier, not only from their scientific work but also from the way they think about music and about what our computational models can, and cannot, capture. Some of the questions developed in this paper would probably not have occurred to me in the same form without these conversations. I am sincerely grateful to both of them.

\medskip

\noindent\textbf{Funding.}
The author is partially supported by the Basque Government through the
BERC 2026--2029 program, by the research project PID2024-155550NB-I00
funded by MICIU/AEI/\allowbreak10.13039/501100011033 and FEDER/EU, and by the
project IT1875-26 funded by the Basque Government.
%%%%%%%%%%%%%%%%%%%%%%%%%%%%%%%%%%%%%%%%%%%%%%%%%%
%%%%%%%%%%%%%%%%%%%%%%%%%%%%%%%%%%%%%%%%%%%%%%%%%%
%%%%%%%%%%%%%%%%%%%%%%%%%%%%%%%%%%%%%%%%%%%%%%%%%%
%%%%%%%%%%%%%%%%%%%%%%%%%%%%%%%%%%%%%%%%%%%%%%%%%%

\end{document}